# ON THE LIMITING DISTRIBUTIONS OF MULTIVARIATE DEPTH-BASED RANK SUM STATISTICS AND RELATED TESTS

BY YIJUN ZUO[1] AND XUMING HE[2]

*Michigan State University and University of Illinois*

A depth-based rank sum statistic for multivariate data introduced by Liu and Singh [*J. Amer. Statist. Assoc.* **88** (1993) 252–260] as an extension of the Wilcoxon rank sum statistic for univariate data has been used in multivariate rank tests in quality control and in experimental studies. Those applications, however, are based on a conjectured limiting distribution, provided by Liu and Singh [*J. Amer. Statist. Assoc.* **88** (1993) 252–260]. The present paper proves the conjecture under general regularity conditions and, therefore, validates various applications of the rank sum statistic in the literature. The paper also shows that the corresponding rank sum tests can be more powerful than Hotelling's $T^2$ test and some commonly used multivariate rank tests in detecting location-scale changes in multivariate distributions.

**1. Introduction.** The key idea of data depth is to provide a center-outward ordering of multivariate observations. Points deep inside a data cloud are assigned high depths, while those on the outskirts are assigned lower depths. The depth of a point decreases when the point moves away from the center of the data cloud. Applications of depth-induced ordering are numerous. For example, Liu and Singh [12] generalized, via data depth, the Wilcoxon rank sum statistic to the multivariate setting. Earlier generalizations of the statistic are due to, for example, Puri and Sen [17], Brown and Hettmansperger [2] and Randles and Peters [19]. More recent ones include Choi and Marden [3], Hettmansperger, Möttönen and Oja [7] and Topchii, Tyurin and Oja [23]. A special version of the Liu–Singh depth-based rank

---

Received July 2003; revised July 2005.

[1]Supported in part by NSF Grants DMS-00-71976 and DMS-01-34628.

[2]Supported in part by NSF Grant DMS-01-02411, NSA Grant H98230-04-1-0034 and NIH Grant R01 DC005603.

---

sum statistic (with a reference sample) inherits the distribution-free property of the Wilcoxon rank sum statistic. The statistic discussed in this paper, like most other generalizations, is only asymptotically distribution-free under the null hypothesis. For its applications in quality control and experimental studies to detect quality deterioration and treatment effects, we refer to [10, 11, 12]. These applications relied on a conjectured limiting distribution, provided by Liu and Singh [12], of the depth-based rank sum statistic. Rousson [20] made an attempt to prove the conjecture, but did not handle the differentiability of the depth functionals for a rigorous treatment. The first objective of the present paper is to fill this mathematical gap by providing regularity conditions for the limiting distribution to hold and by verifying those conditions for some commonly used depth functions. Empirical process theory and, in particular, a generalized Dvoretzk–Kiefer–Wolfowitz theorem in the multivariate setting, turns out to be very useful here.

Our second objective is to investigate the power behavior of the test based on the Liu–Singh rank sum statistic. The test can outperform Hotelling's $T^2$ test and some other existing multivariate tests in detecting location-scale changes for a wide range of distributions. In particular, it is very powerful for detecting scale changes in the alternative, for which Hotelling's $T^2$ test is not even consistent.

Section 2 presents the Liu–Singh depth-based rank sum statistic and an asymptotic normality theorem. Technical proofs of the main theorem and auxiliary lemmas are given in Section 3. The theorem is applied to several commonly used depth functions in Section 4, Section 5 is devoted to a study of the power properties of the rank sum test. Concluding remarks in Section 6 end the paper.

## 2. Liu–Singh statistic and its limiting distribution.
Let $X \sim F$ and $Y \sim G$ be two independent random variables in $\mathbb{R}^d$. Let $D(y; H)$ be a depth function of a given distribution $H$ in $\mathbb{R}^d$ evaluated at point $y$. Lin and Singh [12] introduced $R(y; F) = P_F(X : D(X; F) \leq D(y; F))$ to measure the relative outlyingness of $y$ with respect to $F$ and defined a *quality index*

$$
\begin{aligned}
Q(F, G) &:= \int R(y; F) \, dG(y) \\
&= P\{D(X; F) \leq D(Y; F) \mid X \sim F, Y \sim G\}.
\end{aligned}
\tag{2.1}
$$

Since $R(y; F)$ is the fraction of the $F$ population that is "not as deep" as the point $y$, $Q(F, G)$ is the average fraction over all $y \in G$. As pointed out by Proposition 3.1 of Liu and Singh [12], $R(Y; F) \sim U[0, 1]$ and, consequently, $Q(F, G) = 1/2$ when $Y \sim G = F$ and $D(X; F)$ has a continuous distribution. Thus, the index $Q(F, G)$ can be used to detect a treatment effect or quality deterioration. The Liu–Singh *depth-based rank sum statistic*

$$(2.2) \qquad Q(F_m, G_n) := \int R(y; F_m) \, dG_n(y) = \frac{1}{n} \sum_{j=1}^{n} R(Y_j; F_m)$$

is a two-sample estimator of $Q(F, G)$ based on the empirical distributions $F_m$ and $G_n$. Under the null hypothesis $F = G$ (e.g., no treatment effect or quality deterioration), Liu and Singh [12] proved in *one* dimension $d = 1$ that

$$(2.3) \qquad ((1/m + 1/n)/12)^{-1/2} (Q(F_m, G_n) - 1/2) \xrightarrow{d} N(0, 1)$$

and in *higher* dimensions, they proved the same for the Mahalanobis depth under the existence of the fourth moments, conjecturing that the same limiting distribution holds for general depth functions and in the general multivariate setting. In the next section, we prove this conjecture under some regularity conditions and generalize the result to the case $F \neq G$ in order to perform a power study.

We first list assumptions that are needed for the main result. They will be verified in this and later sections for some commonly used depth functions. Assume, without loss of generality, that $m \leq n$ hereafter. Let $F_m$ be the empirical version of $F$ and $D(\cdot; \cdot)$ be a given depth function with $0 \leq D(x; H) \leq 1$ for any point $x$ and distribution $H$ in $\mathbb{R}^d$.

(A1) $P\{y_1 \leq D(Y; F) \leq y_2\} \leq C|y_2 - y_1|$ for some $C$ and any $y_1, y_2 \in [0, 1]$.
(A2) $\sup_{x \in \mathbb{R}^d} |D(x; F_m) - D(x; F)| = o(1)$, almost surely, as $m \to \infty$.
(A3) $E(\sup_{x \in \mathbb{R}^d} |D(x; F_m) - D(x; F)|) = O(m^{-1/2})$.
(A4) $E(\sum_i p_{iX}(F_m) p_{iY}(F_m)) = o(m^{-1/2})$ if there exist $c_i$ such that $p_{iX}(F_m) > 0$ and $p_{iY}(F_m) > 0$ for $p_{iZ}(F_m) := P(D(Z; F_m) = c_i \mid F_m)$, $i = 1, 2, \ldots$.

Assumption (A1) is the Lipschitz continuity of the distribution of $D(Y; F)$ and can be extended to a more general case with $|x_2 - x_1|$ replaced by $|x_2 - x_1|^\alpha$ for some $\alpha > 0$, if (A3) is also replaced by $E(\sup_{x \in \mathbb{R}^d} |D(x; F_m) - D(x; F)|)^\alpha = O(m^{-\alpha/2})$. The following main result of the paper still holds true:

THEOREM 1. *Let $X \sim F$ and $Y \sim G$ be independent and $X_1, \ldots, X_m$ and $Y_1, \ldots, Y_n$ be independent samples from $F$ and $G$, respectively. Under (A1)–(A4),*

$$(\sigma_{GF}^2/m + \sigma_{FG}^2/n)^{-1/2} (Q(F_m, G_n) - Q(F, G)) \xrightarrow{d} N(0, 1), \qquad \text{as } m \to \infty,$$

*where*

$$\sigma_{FG}^2 = \int P^2(D(X; F) \leq D(y, F)) \, dG(y) - Q^2(F, G),$$

$$\sigma_{GF}^2 = \int P^2(D(x; F) \leq D(Y, F)) \, dF(x) - Q^2(F, G).$$

Assumption (A2) in the theorem is satisfied by most depth functions such as the Mahalanobis, projection, simplicial and halfspace depth functions; see [9, 15, 25] for related discussions. Assumptions (A3)–(A4) also hold true for many of the commonly used depth functions. Verifications can be technically challenging and are deferred to Section 4.

REMARK 1. Under the null hypothesis $F = G$, it is readily seen that $Q(F, G) = 1/2$ and $\sigma_{GF}^2 = \sigma_{FG}^2 = 1/12$ in the theorem.

REMARK 2. Note that (A1)–(A4) and, consequently, the theorem hold true for not only common depth functions that induce a center-outward ordering in $\mathbb{R}^d$, but also other functions that can induce a general (not necessarily center-outward) ordering in $\mathbb{R}^d$. For example, if we define a function $D(x, F) = F(x)$ in $\mathbb{R}^1$, then the corresponding Liu–Singh statistic is equivalent to the Wilcoxon rank sum statistic.

**3. Proofs of the main result and auxiliary lemmas.** To prove the main theorem, we need the following auxiliary lemmas. Some proofs are skipped. For the sake of convenience, we write, for any distribution functions $H$, $F_1$ and $F_2$ in $\mathbb{R}^d$, points $x$ and $y$ in $\mathbb{R}^d$ and a given (affine invariant) depth function $D(\cdot\,;\cdot)$,

$$I(x, y, H) = I\{D(x; H) \leq D(y; H)\},$$

$$I(x, y, F_1, F_2) = I(x, y, F_1) - I(x, y, F_2).$$

LEMMA 1. *Let $F_m$ and $G_n$ be the empirical distributions based on independent samples of sizes $m$ and $n$ from distributions $F$ and $G$, respectively. Then*

   (i)  $\iint I(x, y, F)\, d(G_n(y) - G(y))\, d(F_m(x) - F(x)) = O_p(1/\sqrt{mn})$,
   (ii)  $\iint I(x, y, F_m, F)\, d(F_m(x) - F(x))\, dG(y) = o_p(1/\sqrt{m})$ *under* (A1)–(A2) *and*
   (iii)  $\iint I(x, y, F_m, F)\, dF_m(x)\, d(G_n - G)(y) = O_p(m^{-1/4}n^{-1/2})$ *under* (A1) *and* (A3).

PROOF. We prove (iii). The proofs of (i)–(ii) are omitted. Let $I_{mn} := \iint I(x, y, F_m, F)\, dF_m(x)\, d(G_n - G)(y)$. Then

$$E(I_{mn})^2 \leq E\left\{\int \left[\int I(x, y, F_m, F)\, d(G_n - G)(y)\right]^2 dF_m(x)\right\}$$

$$= E\left[\int I(X_1, y, F_m, F)\, d(G_n - G)(y)\right]^2$$

$$= E\left[ E\left\{ \left( \frac{1}{n} \sum_{j=1}^{n} \left( I(X_1, Y_i, F_m, F) \right. \right.\right.\right.$$
$$\left.\left.\left.\left. - E_Y I(X_1, Y, F_m, F) \right) \right)^2 \Big| X_1, \ldots, X_m \right\} \right]$$

$$\leq \frac{1}{n} E[E_Y\{(I(X_1, Y_1, F_m, F))^2 | X_1, \ldots, X_m\}]$$

$$= \frac{1}{n} E[E_Y\{|I(X_1, Y_1, F_m, F)||X_1, \ldots, X_m\}].$$

One can verify that

$$|I(x, y, F_m, F)| \leq I\left( |D(x; F) - D(y; F)| \leq 2 \sup_{x \in \mathbb{R}^d} |D(x; F_m) - D(x; F)| \right).$$

By (A1) and (A3), we have

$$E(I_{mn})^2 \leq \frac{4C}{n} E\left( \sup_{x \in \mathbb{R}^d} |D(x; F_m) - D(x; F)| \right) = O(1/(m^{1/2}n)).$$

The desired result then follows from Markov's inequality. □

LEMMA 2. *Assume that $X \sim F$ and $Y \sim G$ are independent. Then under* (A4), *we have $\iint I(D(x; F_m) = D(y; F_m))\, dF(x)\, dG(y) = o(m^{-1/2})$.*

PROOF. Let $I(F_m) = \iint I(D(x; F_m) = D(y; F_m))\, dF(x)\, dG(y)$. Conditionally on $X_1, \ldots, X_m$ (or equivalently on $F_m$), we have

$$I(F_m) = \int_{\{y : P(D(X; F_m) = D(y; F_m)|F_m) > 0\}} P(D(X; F_m) = D(y; F_m) \mid F_m)\, dG(y)$$

$$= \sum_i \int_{\{y : P(D(X; F_m) = D(y; F_m) = c_i|F_m) > 0\}} P(D(X; F_m) = c_i \mid F_m)\, dG(y)$$

$$= \sum_i P(D(X; F_m) = c_i \mid F_m) P(D(Y; F_m) = c_i \mid F_m)$$

$$= \sum_i p_{iX}(F_m) p_{iY}(F_m),$$

where $0 \leq c_i \leq 1$ such that $P(D(X; F_m) = c_i \mid F_m) = P(D(Y; F_m) = c_i \mid F_m) > 0$. (Note that there are at most countably many such $c_i$'s.) Taking expectation with respect to $X_1, \ldots, X_m$, the desired result follows immediately from (A4). □

LEMMA 3. *Let $X \sim F$ and $Y \sim G$ be independent and let $X_1, \ldots, X_m$ and $Y_1, \ldots, Y_m$ be independent samples from $F$ and $G$, respectively. Under*

(A1)–(A4),

$$Q(F_m, G_n) - Q(F, G_n) = \iint I(x, y, F) \, dG(y) \, d(F_m(x) - F(x)) + o_p(m^{-1/2})$$

*and, consequently,* $\sqrt{m}(Q(F_m, G_n) - Q(F, G_n)) \overset{d}{\longrightarrow} N(0, \sigma_{GF}^2)$.

PROOF. It suffices to consider the case $F = G$. First, we observe that

$$Q(F_m, G_n) - Q(F, G_n) = \int R(y; F_m) \, dG_n(y) - \int R(y; F) \, dG_n(y)$$

$$= \iint I(x, y, F_m) \, dF_m(x) \, dG_n(y)$$

$$- \iint I(x, y, F) \, dF(x) \, dG_n(y)$$

$$= \iint [I(x, y, F_m) - I(x, y, F)] \, dF_m(x) \, dG_n(y)$$

$$+ \iint I(x, y, F) \, d(G_n(y) - G(y)) \, d(F_m(x) - F(x))$$

$$+ \iint I(x, y, F) \, dG(y) \, d(F_m(x) - F(x)).$$

We shall call the last three terms $I_{mn1}$, $I_{mn2}$ and $I_{m3}$, respectively. From Lemma 1, it follows immediately that $\sqrt{m} \, I_{mn2} = o_p(1)$. By a standard central limit theorem, we have

$$(3.4) \qquad\qquad \sqrt{m} \, I_{m3} \overset{d}{\longrightarrow} N(0, \sigma_{GF}^2).$$

We now show that $\sqrt{m} \, I_{mn1} = o_p(1)$. Observe that

$$I_{mn1} = \iint I(x, y, F_m, F) \, dF_m(x) \, d(G_n - G)(y)$$

$$+ \iint I(x, y, F_m, F) \, dF_m(x) \, dG(y)$$

$$= \iint I(x, y, F_m, F) \, dF(x) \, dG(y) + o_p(1/\sqrt{m}),$$

by Lemma 1 and the given condition. It is readily seen that

$$\iint I(x, y, F_m, F) \, dF(x) \, dG(y)$$

$$= \iint I(x, y, F_m) \, dF(x) \, dG(y) - \iint I(x, y, F) \, dF(x) \, dG(y)$$

$$= \frac{1}{2} \iint [I(D(x, F_m) \leq D(y; F_m))$$

$$+ I(D(x, F_m) \geq D(y; F_m))] \, dF(x) \, dG(y) - \frac{1}{2}$$

$$= \frac{1}{2} \iint I(D(x, F_m) = D(y; F_m)) \, dF(x) \, dG(y) = o(m^{-1/2}),$$

by Lemma 2. The desired result follows immediately. $\square$

Proof of Theorem 1.   By Lemma 3, we have

$$Q(F_m, G_n) - Q(F, G)$$
$$= (Q(F_m, G_n) - Q(F, G_n)) + (Q(F, G_n) - Q(F, G))$$
$$= \iint I(x, y, F) \, dG(y) \, d(F_m(x) - F(x))$$
$$\quad + \iint I(x, y, F) \, dF(x) \, d(G_n(y) - G(y)) + o_p(m^{-1/2}).$$

The independence of $F_m$ and $G_n$ and the central limit theorem then give the result. $\square$

## 4. Applications and examples.

This section verifies (A3)–(A4) [and (A2)] for several common depth functions. Mahalanobis, halfspace and projection depth functions are selected for illustration. The findings here and in Section 2 ensure the validity of Theorem 1 for these depth functions.

Example 1 [Mahalanobis depth (MHD)].   The depth of a point $x$ is defined as

$$MHD(x; F) = 1/(1 + (x - \mu(F))'\Sigma^{-1}(F)(x - \mu(F))), \qquad x \in \mathbb{R}^d,$$

where $\mu(F)$ and $\Sigma(F)$ are location and covariance measures of a given distribution $F$; see [12, 27]. Clearly, both $MHD(x; F)$ and $MHD(x; F_m)$ vanish at infinity as $\|x\| \to \infty$, where $F_m$ is the empirical version of $F$ based on $X_1, \ldots, X_m$ and $\mu(F_m)$ and $\Sigma(F_m)$ are strongly consistent estimators of $\mu(F)$ and $\Sigma(F)$, respectively. Hence,

$$\sup_{x \in \mathbb{R}^d} |MHD(x; F_m) - MHD(x; F)| = |MHD(x_m; F_m) - MHD(x_m; F)|,$$

by the continuity of $MHD(x; F)$ and $MHD(x; F_m)$ in $x$ for some $x_m = x(F_m, F) \in \mathbb{R}^d$ such that $\|x_m\| \le M < \infty$ for some $M > 0$ and all large $m$. Write, for simplicity, $\mu$ and $\Sigma$ for $\mu(F)$ and $\Sigma(F)$ and $\mu_m$ and $\Sigma_m$ for $\mu(F_m)$ and $\Sigma(F_m)$, respectively. Then

$$|MHD(x_m; F_m) - MHD(x_m; F)|$$
$$= \frac{|(\mu_m - \mu)'\Sigma_m^{-1}(\mu_m + \mu - 2x_m) + (x_m - \mu)'(\Sigma_m^{-1} - \Sigma^{-1})(x_m - \mu)|}{(1 + \|\Sigma_m^{-1/2}(x_m - \mu_m)\|^2)(1 + \|\Sigma^{-1/2}(x_m - \mu)\|^2)}.$$

This, in conjunction with the strong consistency of $\mu_m$ and $\Sigma_m$, yields (A2).

Hölder's inequality and expectations of quadratic forms (page 13 of [21]) yield (A3) if conditions (i) and (ii) below are met. (A4) holds trivially if (iii) holds.

(i) $\mu_m$ and $\Sigma_m$ are strongly consistent estimators of $\mu$ and $\Sigma$, respectively;

(ii) $E(\mu_m - \mu)_i = O(m^{-1/2})$, $E(\Sigma_m^{-1} - \Sigma^{-1})_{jk} = O(m^{-1/2})$, $1 \le i, j, k \le d$, where the subscripts $i$ and $jk$ denote the elements of a vector and a matrix, respectively;

(iii) The probability mass of $X$ at any ellipsoid is 0.

COROLLARY 1. *Assume that conditions* (i), (ii) *and* (iii) *hold and the distribution of $MHD(Y; F)$ is Lipschitz continuous. Then Theorem* 1 *holds for MHD.*

EXAMPLE 2 [Halfspace depth (HD)]. Tukey [24] suggested this depth as

$$HD(x; F) = \inf\{P(H_x) : H_x \text{ closed halfspace with } x \text{ on its boundary}\},$$
$$x \in \mathbb{R}^d,$$

where $P$ is the probability measure corresponding to $F$. (A2) follows immediately (see, e.g., pages 1816–1817 of [5]). Let $\mathcal{H}$ be the set of all closed halfspaces and $P_m$ be the empirical probability measure of $P$. Define

$$D_m(\mathcal{H}) := m^{1/2}\|P_m - P\|_{\mathcal{H}} := \sup_{H \in \mathcal{H}} m^{1/2}|P_m(H) - P(H)|.$$

Note that $\mathcal{H}$ is a permissible class of sets with polynomial discrimination (see Section II.4 of [18] for definitions and arguments). Let $S(\mathcal{H})$ be the degree of the corresponding polynomial. Then, by a generalized Dvoretzky–Kiefer–Wolfowitz theorem (see [1, 13, 14] and Section 6.5 of [6]), we have, for any $\varepsilon > 0$, that $P(D_m(\mathcal{H}) > M) \le Ke^{-(2-\varepsilon)M^2}$ for some sufficiently large constant $K = K(\varepsilon, S(\mathcal{H}))$. This immediately yields (A3).

To verify (A4), we consider the case $F = G$ for simplicity. We first note that $HD(X; F_m)$ for given $F_m$ is discrete and can take at most $O(m)$ values $c_i = i/m$ for $i = 0, 1, \ldots, m$. Let $F$ be continuous. We first consider the univariate case. Let

$$A_0 = \mathbb{R}^1 - \bigcap H_m, \qquad A_i = \bigcap H_{m-i+1} - \bigcap H_{m-i},$$

$$A_{k+1} = \bigcap H_{m-k} - \bigcap \varnothing,$$

with $1 \leq i \leq k$ and $k = \lfloor (m-1)/2 \rfloor$, where $H_i$ is any closed half-line containing exactly some $i$ points of $X_1, \ldots, X_m$. It follows that for $0 \leq i \leq k$,

$$P(HD(X; F_m) = c_i \mid F_m)$$
$$= P(A_i) = [F(X_{(i+1)}) - F(X_{(i)})] + [F(X_{(m-i+1)}) - F(X_{(m-i)})]$$
$$P(HD(X; F_m) = c_{k+1} \mid F_m) = P(A_{k+1}) = [F(X_{(m-k)}) - F(X_{(k+1)})],$$

where $-\infty =: X_{(0)} \leq X_{(1)} \leq \cdots \leq X_{(m)} \leq X_{(m+1)} := \infty$ are order statistics.

On the other hand, $X_{(i)}$ and $F^{-1}(U_{(i)})$ are equal in distribution ($\overset{d}{=}$), where $0 =: U_{(0)} \leq U_{(1)} \leq \cdots \leq U_{(m)} \leq U_{(m+1)} := 1$ are order statistics based on a sample from the uniform distribution on $[0,1]$. Let $D_i = F(X_{(i+1)}) - F(X_{(i)})$, $i = 0, \ldots, m$. The $D_i'$s have the same distribution and

$$E(D_i) = \frac{1}{m+1}, \qquad E(D_i^2) = \frac{2}{(m+1)(m+2)},$$

$$E(D_i D_j) = \frac{1}{(m+1)(m+2)}.$$

Hence, for $0 \leq i \leq k$, $E((P(HD(X; F_m) = c_i \mid F_m))^2) = 6/((m+1)(m+2))$ and $E((P(HD(X; F_m) = c_{k+1} \mid F_m))^2) = O(m^{-2})$. Thus, (A4) follows immediately.

Let us now treat the multivariate case. Let $X_1, \ldots, X_m$ be given. Denote by $H_i$ any closed halfspace containing exactly $i$ points of $X_1, \ldots, X_m$. Define the sets

$$A_0 = \mathbb{R}^d - \bigcap H_m, \qquad A_1 = \bigcap H_m - \bigcap H_{m-1} \qquad, \ldots,$$

$$A_{m-k} = \bigcap H_{k+1} - \bigcap H_k, \qquad A_{m-k+1} = \bigcap H_k,$$

with $(m-k+1)/m = \max_{x \in \mathbb{R}^d} HD(x; F_m) \leq 1$. Then it is not difficult to see that

$$HD(x; F_m) = i/m, \qquad \text{for } x \in A_i, \ i = 0, 1, \ldots, m-k, m-k+1.$$

Now let $p_i = P(HD(X; F_m) = c_i \mid F_m)$ with $c_i = i/m$. Then, for any $0 \leq i \leq m-k+1$, $p_i = P(X \in A_i) = P(\bigcap H_{m-i+1}) - P(\bigcap H_{m-i})$ with $H_{m+1} = \mathbb{R}^d$ and $H_{k-1} = \varnothing$. Now, treating $p_i$ as random variables based on the random variables $X_1, \ldots, X_m$, by symmetry and the uniform spacings results used for the univariate case above, we conclude that the $p_i$'s have the same distribution for $i = 0, \ldots, m-k$ and

$$E(p_i) = O(m^{-1}), \quad E(p_i^2) = O(m^{-2}), \qquad i = 0, \ldots, m-k+1.$$

Assumption (A4) follows in a straightforward fashion. Thus, we have COROLLARY 2. *Assume that $F$ is continuous and the distribution of $HD(Y; F)$ is Lipschitz continuous. Then Theorem* 1 *holds true for HD.*

EXAMPLE 3 [Projection depth (PD)]. Stahel [22] and Donoho [4] defined the outlyingness of a point $x \in \mathbb{R}^d$ with respect to $F$ in $\mathbb{R}^d$ as

$$O(x; F) = \sup_{u \in S^{d-1}} |u'x - \mu(F_u)|/\sigma(F_u),$$

where $S^{d-1} = \{u : \|u\| = 1\}$, $\mu(\cdot)$ and $\sigma(\cdot)$ are univariate location and scale estimators such that $\mu(aZ + b) = a\mu(Z) + b$ and $\sigma(aZ + b) = |a|\sigma(Z)$ for any scalars $a, b \in \mathbb{R}^1$ and random variable $Z \in \mathbb{R}^1$ and $u'X \sim F_u$ with $X \sim F$. The projection depth of $x$ with respect to $F$ is then defined as (see [10, 25])

$$PD(x; F) = 1/(1 + O(x; F)).$$

Under the following conditions on $\mu$ and $\sigma$,

(C1)  $\sup_{u \in S^{d-1}} \mu(F_u) < \infty$, $0 < \inf_{u \in S^{d-1}} \sigma(F_u) \leq \sup_{u \in S^{d-1}} \sigma(F_u) < \infty$;
(C2)  $\sup_{u \in S^{d-1}} |\mu(F_{mu}) - \mu(F_u)| = o(1)$, $\sup_{u \in S^{d-1}} |\sigma(F_{mu}) - \sigma(F_u)| = o(1)$, a.s.
(C3)  $E \sup_{\|u\|=1} |\mu(F_{mu}) - \mu(F_u)| = O(m^{-1/2})$, $E \sup_{\|u\|=1} |\sigma(F_{mu}) - \sigma(F_u)| = O(m^{-1/2})$,

where $F_{mu}$ is the empirical distribution based on $u'X_1, \ldots, u'X_m$ and $X_1, \ldots, X_m$ is a sample from $F$, Assumption (A2) holds true by Theorem 2.3 of [25] and (A3) follows from (C3) and the fact that for any $x \in \mathbb{R}^d$ and some constant $C > 0$,

$$|PD(x; F_m) - PD(x; F)|$$
$$\leq \sup_{u \in S^{d-1}} \frac{O(x; F)|\sigma(F_{mu}) - \sigma(F_u)| + |\mu(F_{mu}) - \mu(F_u)|}{(1 + O(x; F_m))(1 + O(x; F))\sigma(F_{mu})}$$
$$\leq C \sup_{u \in S^{d-1}} \{|\sigma(F_{mu}) - \sigma(F_u)| + |\mu(F_{mu}) - \mu(F_u)|\}.$$

(C1)–(C3) is true for general smooth $M$-estimators of $\mu$ and $\sigma$ (see [8]) and rather general distribution functions $F$. If we consider the median (Med) and the median absolute deviation (MAD), then (C3) holds under the following condition:

(C4)  $F_u$ has a continuous density $f_u$ around points $\mu(F_u) + \{0, \pm\sigma(F_u)\}$ such that

$$\inf_{\|u\|=1} f_u(\mu(F_u)) > 0,$$

$$\inf_{\|u\|=1} (f_u(\mu(F_u) + \sigma(F_u)) + f_u(\mu(F_u) - \sigma(F_u))) > 0.$$

To verify this, it suffices to establish just the first part of (C3) for $\mu = \text{Med}$. Observe that

$$
\begin{aligned}
F_u^{-1}(1/2 - \|F_{mu} - F_u\|_\infty) - F_u^{-1}(1/2) &\leq \mu(F_{mu}) - \mu(F_u) \\
&\leq F_u^{-1}(1/2 + \|F_{mu} - F_u\|_\infty) \\
&\quad - F_u^{-1}(1/2)
\end{aligned}
$$

for any $u$ and sufficiently large $m$. Hence,

$$
|\mu(F_{mu}) - \mu(F_u)| \leq 2\|F_{mu} - F_u\|_\infty / \inf_{u \in S^{d-1}} f_u(\mu(F_u)) := C\|F_{mu} - F_u\|_\infty,
$$

by (C4). Clearly, $\mu(F_{mu})$ is continuous in $u$. From (C4), together with and Lemma 5.1 and Theorem 3.3 of [25], it follows that $\mu(F_u)$ is also continuous in $u$. Therefore,

$$
\begin{aligned}
P\Big(\sqrt{m} \sup_{u \in S^{d-1}} |\mu(F_{mu}) - \mu(F_u)| > t\Big) &\leq P(\|F_{mu_0} - F_{u_0}\|_\infty > (t^2/(mC^2))^{1/2}) \\
&\leq 2e^{-2t^2/C^2}, \qquad \text{for any } t > 0,
\end{aligned}
$$

where the unit vector $u_0$ may depend on $m$. Hence, the first part of (C3) follows.

Assumption (A4) holds for $PD$ since $P(PD(X; F_m) = c \mid F_m) = 0$ for most commonly used $(\mu, \sigma)$ and $F$. First, the continuity of $\mu(F_{mu})$ and $\sigma(F_{mu})$ in $u$ gives

$$
P(PD(X; F_m) = c \mid F_m) = P((u_X' X - \mu(F_{mu_X}))/\sigma(F_{mu_X}) = (1-c)/c \mid F_m)
$$

for some unit vector $u_X$ depending on $X$. This probability is 0 for most $F$ and $(\mu, \sigma)$. For example, if $(\mu, \sigma) = (\text{mean, standard deviation})$, then

$$
P(PD(X; F_m) = c \mid F_m) = P(\|S_m^{-1/2}(X - \bar{X}_m)\| = (1-c)/c \mid F_m),
$$

where $S_m = \frac{1}{m-1} \sum_{i=1}^m (X_i - \bar{X}_m)(X_i - \bar{X}_m)'$, which is 0 provided the mass of $F$ on any ellipsoid is 0. Thus,

COROLLARY 3. *Assume that* (C1)–(C3) *hold,* $P(PD(X; F_m) = c \mid F_m) = 0$ *for any* $c \geq 0$ *and* $PD(Y; F)$ *satisfies* (A1). *Then Theorem* 1 *holds for* $PD$.

**5. Power properties of the Liu–Singh multivariate rank sum test.**

*Large sample properties.* A major application of the Liu–Singh multivariate rank-sum statistic is to test the following hypotheses:

$$(5.5) \qquad H_0 : F = G \qquad \text{versus} \qquad H_1 : F \neq G.$$

By Theorem 1, a large sample test based on the Liu–Singh rank-sum statistic $Q(F_m, G_n)$ rejects $H_0$ at (an asymptotic) significance level $\alpha$ when

$$(5.6) \qquad |Q(F_m, G_n) - 1/2| > z_{1-\alpha/2}((1/m + 1/n)/12)^{1/2},$$

where $\Phi(z_r) = r$ for $0 < r < 1$ and normal cumulative distribution function $\Phi(\cdot)$. The test is *affine invariant* and is *distribution-free* in the asymptotic sense under the null hypothesis. Here, we focus on the asymptotic power properties of the test. By Theorem 1, the (asymptotic) power function of the depth-based rank-sum test with an asymptotic significance level $\alpha$ is

$$
(5.7)
\begin{aligned}
\beta_Q(F, G) = {} & 1 - \Phi\left( \frac{1/2 - Q(F, G) + z_{1-\alpha/2}\sqrt{(1/m + 1/n)/12}}{\sqrt{\sigma_{GF}^2/m + \sigma_{FG}^2/n}} \right) \\
& + \Phi\left( \frac{1/2 - Q(F, G) - z_{1-\alpha/2}\sqrt{(1/m + 1/n)/12}}{\sqrt{\sigma_{GF}^2/m + \sigma_{FG}^2/n}} \right).
\end{aligned}
$$

The asymptotic power function indicates that the test is *consistent* for *all* alternative distributions $G$ such that $Q(F, G) \neq 1/2$. Before studying the behavior of $\beta_Q(F, G)$, we shall consider its key component $Q(F, G)$, the so-called *quality index* in [12]. For convenience, consider a normal family and let $d = 2$. Assume, without loss of generality, that $F = N_2((0, 0)', I_2)$ and consider $G = N_2(\mu, \Sigma)$, where $I_2$ is the $2 \times 2$ identity matrix. It can be shown that

$$Q(F, G) = (|S|/|\Sigma|)^{1/2} \exp(-\mu'(\Sigma^{-1} - \Sigma^{-1} S \, \Sigma^{-1})\mu/2)$$

for *any* affine invariant depth functions, where $S = (I_2 + \Sigma^{-1})^{-1}$. In the case $\mu = (u, u)'$ and $\Sigma = \sigma^2 I_2$, write $Q(u, \sigma^2)$ for $Q(F, G)$. Then

$$Q(u, \sigma^2) := Q(F, G) = \exp(-u^2/(1 + \sigma^2))/(1 + \sigma^2).$$

Its behavior is revealed in Figure 1. It increases to its maximum value $(1 + \sigma^2)^{-1}$ [or $\exp(-u^2)$] as $u \to 0$ for a fixed $\sigma^2$ (or as $\sigma^2 \to 0$ for a fixed $u$). When $u = 0$ and $\sigma^2 = 1$, $Q(F, G)$, as expected, is $1/2$, and it is less than $1/2$ when there is a dilution in the distribution ($\sigma^2 > 1$). Note that Liu and Singh [12] also discussed $Q(u, \sigma^2)$. The results here are more accurate than their Table 1.

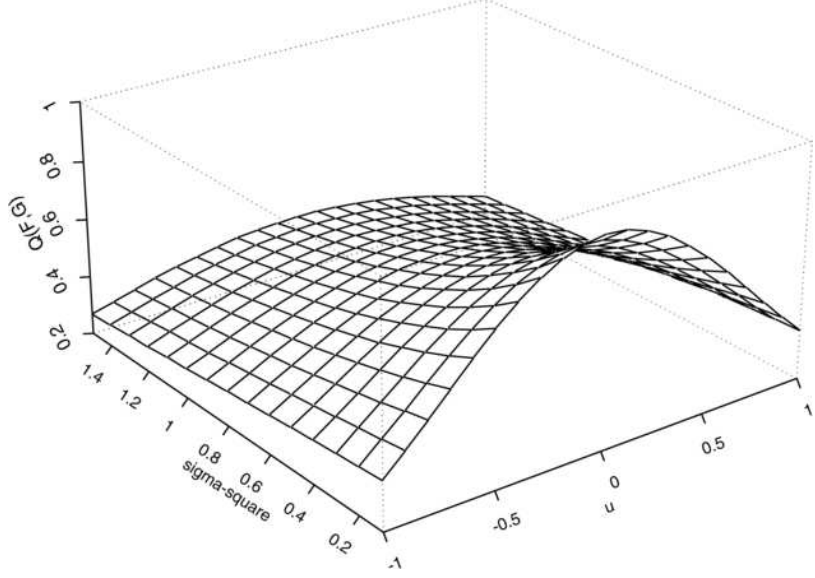

Fig. 1. *The behavior of $Q(F, G)$ with $F = N_2((0, 0)', I_2)$ and $G = N_2((u, u)', \sigma^2 I_2)$.*

A popular large sample test for hypotheses (5.5) is based on Hotelling's $T^2$ statistic that rejects $H_0$ if

$$(5.8) \qquad (\bar{X}_m - \bar{Y}_n)'((1/m + 1/n)S_{pooled})^{-1}(\bar{X}_m - \bar{Y}_n) > \chi^2_{1-\alpha}(d),$$

where $S_{pooled} = ((m-1)S_X + (n-1)S_Y)/(m+n-2)$, $\bar{X}_m$, $\bar{Y}_n$, $S_X$ and $S_Y$ are sample means and covariance matrices and $\chi^2_r(d)$ is the $r$th quantile of the chi-square distribution with $d$ degrees of freedom. The power function of the test is

$$
\begin{aligned}
(5.9) \qquad \beta_{T^2}(F, G) = P(&(\bar{X}_m - \bar{Y}_n)'((1/m + 1/n)S_{pooled})^{-1}(\bar{X}_m - \bar{Y}_n) \\
&> \chi^2_{1-\alpha}(d)).
\end{aligned}
$$

We also consider a multivariate rank-sum test based on the Oja objective function in [7]. The Oja test statistic, $O$, has the following null distribution with $N = m + n$ and $\lambda = n/N$:

$$(5.10) \qquad O := (N\lambda(1-\lambda))^{-1} T'_N B_N^{-1} T_N \overset{d}{\longrightarrow} \chi^2(d).$$

Here $T_N = \sum_{k=1}^N a_k R_N(z_k)$, $R_N(z) = \frac{d!(N-d)!}{N!} \sum_{p \in P} S_p(z) n_p$, $a_k = (1-\lambda)I(k > m) - \lambda I(k < m)$, $z_k \in \{X_1, \dots, X_m, Y_1, \dots, Y_n\}$, $B_N = \frac{1}{N-1} \sum_k R_N(z_k) R'_N(z_k)$, $P = \{p = (i_1, \dots, i_d) : 1 \leq i_1 < \dots < i_d \leq N\}$, $S_p(z) = \text{sign}(n_{0p} + z' n_p)$ and

$$\det \begin{pmatrix} 1 & 1 & \dots & 1 & 1 \\ z_{i_1} & z_{i_2} & \dots & z_{i_d} & z \end{pmatrix} = n_{0p} + z' n_p,$$

where $n_{0p}$ and $n_{jp}, j = 1, \ldots, d$, are the cofactors according to the column $(1, z')'$. The power function of this rank test with an asymptotic significance level $\alpha$ is

(5.11)        $\beta_O(F, G) = P((N\lambda(1-\lambda))^{-1}T_N'B_N^{-1}T_N > \chi_{1-\alpha}^2(d))$.

The asymptotic relative efficiency (ARE) of this test in Pitman's sense is discussed in the literature; see, for example, [16]. At the bivariate normal model, it is 0.937 relative to $T^2$.

In the following, we study the behavior of $\beta_Q$, $\beta_O$ and $\beta_{T^2}$. To facilitate our discussion, we assume that $\alpha = 0.05$, $m = n$, $d = 2$ and that $G$ is normal or mixed (contaminated) normal, shrinking to the null distribution $F = N_2((0,0)', I_2)$. Note that the asymptotic power of the depth-based rank-sum test, hereafter called the $Q$ test, is invariant with respect to the choice of the depth function.

For pure location shift models $Y \sim G = N_2((u, u)', I_2)$, Hotelling's $T^2$ based test, hereafter called $T^2$ is the most powerful, followed by the Oja rank test, to be called the $O$ test, and then followed by the $Q$ test. All of these tests are consistent at any fixed alternative. Furthermore, we note that when the dimension $d$ gets larger, the asymptotic powers of these tests move closer.

On the other hand, for pure scale change models $G = N_2((0,0)', \sigma^2 I_2)$, the $Q$ test is much more powerful than the other test. In fact, for these models, the $T^2$ test has trivial asymptotic power $\alpha$ at all alternatives. Figure 2, a plot of the power functions $\beta_{T^2}$, $\beta_O$ and $\beta_Q$, clearly reveals the superiority of the $Q$ test. The $O$ test performs just slightly better than $T^2$.

In the following, we consider a location shift with contamination a scale change with contamination and a simultaneous location and scale change as alternatives. The amount of contamination $\varepsilon$ is set to be 10%. The asymptotic power calculations for $T^2$ and $Q$ are based on the limiting distributions of the test statistics under the alternatives. Since the limiting distribution is not available for the $O$ test (except for pure location shift models), we use Monte Carlo to estimate the powers.

For $G = (1-\varepsilon)N_2((u, u)', I_2) + \varepsilon N_2((0,0)', (1 + 10u\sigma^2)I_2)$, the contaminated location shift models with $u \geq 0$, the (*asymptotic*) power function $\beta_{T^2}(F, G)$ is $P(Z_{2a} \geq \chi_{0.95}^2(2))$, where $Z_{2a}$ has a noncentral chi-square distribution with two degrees of freedom and noncentrality parameter $n(1 - \varepsilon)^2 u^2 / (1 + 5\varepsilon u\sigma^2 + \varepsilon(1 - \varepsilon)u^2)$. Since the derivation of this result is quite tedious, we omit the details. Comparisons of $\beta_{T^2}$, $\beta_O$ and $\beta_Q$ are listed in Table 1, which clearly reveals that $T^2$ becomes less powerful than $Q$ when a pure location shift model is 10% contaminated. For large $\mu$, $O$ is more powerful than $Q$ since the underlying model is mainly a location shift.

For $G = 0.9N_2((0,0)', \sigma^2 I_2) + 0.1N_2((u, u)', I_2)$, the contaminated scale change models with $\sigma = u + 1 \geq 1$, the (*asymptotic*) power function $\beta_{T^2}$ is equal to $P(Z_{2b} \geq \chi^2_{0.95}(2))$, where $Z_{2b}$ has the noncentral chi-square distribution with two degrees of freedom and noncentrality parameter $2n\varepsilon^2 u^2/(1 + \varepsilon + (1-\varepsilon)\sigma^2 + 2\varepsilon(1-\varepsilon)u^2)$. Table 2 reveals the superiority of $Q$ in detecting scale changes over $T^2$ and $O$, even when the model has a 10% contamination.

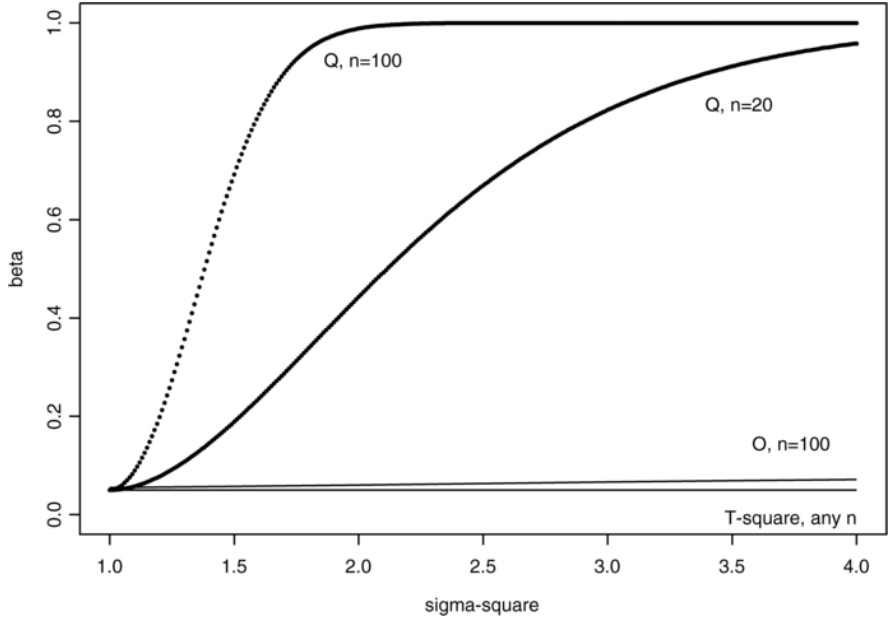

Fig. 2. $\beta_{T^2}(F, G)$ and $\beta_Q(F, G)$ with $F = N_2((0,0)', I_2)$, $G = N_2((0,0)', \sigma^2 I_2)$.

TABLE 1
The (asymptotic) power of tests based on $T^2$, $O$ and $Q$.

| $u$ | | 0.0 | 0.15 | 0.20 | 0.25 | 0.30 | 0.35 |
|---|---|---|---|---|---|---|---|
| | | $G = 0.9N_2((u,u)', I_2) + 0.1N_2((0,0)', (1 + 10u\sigma^2)I_2)$, $\quad \sigma = 4$ | | | | | |
| $n = 100$ | $\beta_{T^2}$ | 0.050 | 0.117 | 0.155 | 0.196 | 0.239 | 0.284 |
| | $\beta_Q$ | 0.051 | 0.245 | 0.286 | 0.307 | 0.381 | 0.443 |
| | $\beta_O$ | 0.046 | 0.157 | 0.296 | 0.423 | 0.558 | 0.687 |
| $n = 200$ | $\beta_{T^2}$ | 0.050 | 0.193 | 0.273 | 0.357 | 0.441 | 0.521 |
| | $\beta_Q$ | 0.051 | 0.430 | 0.508 | 0.549 | 0.659 | 0.746 |
| | $\beta_O$ | 0.056 | 0.342 | 0.546 | 0.712 | 0.881 | 0.941 |

For $G = N_2((u,u)', \sigma^2 I_2)$, the simultaneous location and scale change models with $\sigma = u + 1 \geq 1$, the (asymptotic) power function $\beta_{T^2}$ is $P(Z_{2c} \geq \chi^2_{0.95}(2))$, where $Z_{2c}$ has the noncentral chi-square distribution with two degrees of freedom and noncentrality parameter $2nu^2/(1 + \sigma^2)$. Table 3 reveals that $Q$ can be more powerful than $T^2$ and $O$ when there are simultaneous location and scale changes. Here, we selected $(\sigma - 1)/u = 1$. Our empirical evidence indicates that the superiority of $Q$ holds provided that $(\sigma - 1)/u$ is close to or greater than 1, that is, as long as the change in scale is not much less than that in location. Also note that in this model, $T^2$ is more powerful than $O$.

*Small sample properties.* To check the small sample power behavior of $Q$, we now examine the empirical behavior of the test based on $Q(F_m, G_n)$ and compare it with those of $T^2$ and $O$. We focus on the relative frequencies of rejecting $H_0$ of (5.5) at $\alpha = 0.05$ based on the tests (5.6), (5.8) and (5.10) and 1000 samples from $F$ and $G$ at the sample size $m = n = 25$. The projection depth with $(\mu, \sigma) = $ (Med, MAD) is selected in our simulation studies and some results are given in Table 4. Again, we skip the pure location shift and scale change models, in which cases, $T^2$ and $Q$ perform best, respectively. Our Monte Carlo studies confirm the validity of the (asymptotic) power properties of $Q$ at small samples.

TABLE 2
*The (asymptotic) power of tests based on $T^2$, $O$ and $Q$.*

| $\sigma^2$ | | **1.0** | **1.2** | **1.4** | **1.6** | **1.8** | **2.0** |
|---|---|---|---|---|---|---|---|
| | | $G = 0.9 N_2((0,0)', \sigma^2 I_2) + 0.1 N_2((u,u)', I_2),$ | | | | $u = \sigma - 1$ | |
| $n = 100$ | $\beta_{T^2}$ | 0.050 | 0.050 | 0.052 | 0.054 | 0.056 | 0.059 |
| | $\beta_Q$ | 0.051 | 0.181 | 0.430 | 0.734 | 0.891 | 0.963 |
| | $\beta_O$ | 0.048 | 0.054 | 0.057 | 0.064 | 0.068 | 0.070 |
| $n = 200$ | $\beta_{T^2}$ | 0.050 | 0.051 | 0.054 | 0.058 | 0.063 | 0.068 |
| | $\beta_Q$ | 0.051 | 0.299 | 0.740 | 0.950 | 0.994 | 1.000 |
| | $\beta_O$ | 0.052 | 0.059 | 0.063 | 0.085 | 0.112 | 0.139 |

TABLE 3
*The (asymptotic) power of tests based on $T^2$, $Q$ and $O$.*

| $u$ | | **0.0** | **0.15** | **0.20** | **0.25** | **0.30** | **0.35** |
|---|---|---|---|---|---|---|---|
| | | $G = N_2((u,u)', \sigma^2 I_2),$ | | $\sigma = u + 1$ | | | |
| $n = 100$ | $\beta_{T^2}$ | 0.050 | 0.219 | 0.348 | 0.493 | 0.634 | 0.755 |
| | $\beta_Q$ | 0.051 | 0.437 | 0.662 | 0.839 | 0.941 | 0.983 |
| | $\beta_O$ | 0.046 | 0.218 | 0.324 | 0.430 | 0.573 | 0.708 |
| $n = 200$ | $\beta_{T^2}$ | 0.050 | 0.404 | 0.625 | 0.805 | 0.916 | 0.970 |
| | $\beta_Q$ | 0.049 | 0.725 | 0.922 | 0.987 | 0.999 | 1.000 |
| | $\beta_O$ | 0.056 | 0.357 | 0.569 | 0.755 | 0.882 | 0.944 |

TABLE 4
*Observed relative frequency of rejecting $H_0$.*

| $u$ | | **0.0** | **0.15** | **0.20** | **0.25** | **0.30** | **0.35** |
|---|---|---|---|---|---|---|---|
| | | $G = 0.9N_2((u,u)', I_2) + 0.1N_2((0,0)', (1+10u\sigma^2)I_2)$, | | | $\sigma = 4$ | | |
| $n = 25$ | $\beta_{T^2}$ | 0.058 | 0.083 | 0.108 | 0.142 | 0.151 | 0.189 |
| | $\beta_Q$ | 0.057 | 0.154 | 0.156 | 0.170 | 0.203 | 0.216 |
| | $\beta_O$ | 0.047 | 0.084 | 0.116 | 0.152 | 0.201 | 0.254 |
| $\sigma^2$ | | **1.0** | **1.2** | **1.4** | **1.6** | **1.8** | **2.0** |
| | | $G = 0.9N_2((0,0)', \sigma^2 I_2) + 0.1N_2((u,u)', I_2)$, | | | $u = \sigma - 1$ | | |
| $n = 25$ | $\beta_{T^2}$ | 0.059 | 0.063 | 0.059 | 0.073 | 0.061 | 0.067 |
| | $\beta_Q$ | 0.063 | 0.145 | 0.243 | 0.377 | 0.469 | 0.581 |
| | $\beta_O$ | 0.051 | 0.058 | 0.041 | 0.053 | 0.043 | 0.055 |
| $u$ | | **0.0** | **0.15** | **0.20** | **0.25** | **0.30** | **0.35** |
| | | $G = N_2((u,u)', \sigma^2 I_2)$, | | | $\sigma = u + 1$ | | |
| $n = 25$ | $\beta_{T^2}$ | 0.069 | 0.113 | 0.147 | 0.183 | 0.220 | 0.269 |
| | $\beta_Q$ | 0.060 | 0.245 | 0.324 | 0.418 | 0.498 | 0.587 |
| | $\beta_O$ | 0.044 | 0.082 | 0.089 | 0.127 | 0.197 | 0.221 |

**6. Concluding remarks.** This paper proves the conjectured limiting distribution of the Liu–Singh multivariate rank-sum statistic under some regularity conditions which are verified for several commonly used depth functions. The asymptotic results in the paper are established for *general* depth structures and for *general* distributions $F$ and $G$. The $Q$ test requires neither the existence of a covariance matrix nor the symmetry of $F$ and $G$. This is not always the case for Hotelling's $T^2$ test and other multivariate generalizations of Wilcoxon's rank-sum test.

The paper also studies the power behavior of the rank-sum test both asymptotically and empirically. Although the discussion focuses on the normal and mixed normal models and $d = 2$, what we learned from these investigations is typical for $d > 2$ and for many non-Gaussian models. Our investigations also indicate that the conclusions drawn from our two-sample problems are valid for one-sample problems.

The Liu–Singh rank-sum statistic plays an important role in detecting scale changes similarly to the role played by Hotelling's $T^2$ in detecting location shifts of distributions. When there is a scale change in $F$, the depths of almost all points $y$ from $G$ decrease or increase together and, consequently, $Q(F, G)$ is very sensitive to the change. This explains why $Q$ is so powerful in detecting small scale changes. On the other hand, when there is a small shift in location, the depths of some points $y$ from $G$ increase, whereas those of the others decrease and, consequently, $Q(F, G)$ will not be so sensitive to a small shift in location. Unlike the $T^2$ test for scale change alternatives, the $Q$ test *is* consistent for location shift alternatives, nevertheless.

Finally, we briefly address the computing issue. Hotelling's $T^2$ statistic is clearly the easiest to compute. The computational complexity for the $O$ test is $O(n^{d+1})$, while the complexity for the $Q$ test based on the projection depth (PD) is $O(n^{d+2})$. Indeed, the projection depth can be computed exactly by considering $O(n^d)$ directions that are perpendicular to a hyperplane determined by $d$ data points; see [26] for a related discussion. The exact computation is of course time-consuming. In our simulation study, we employed approximate algorithms (see http://www.stt.msu.edu/~zuo/table4.txt), which consider a large number of directions in computing $Q$ and $O$.

**Acknowledgments.** The authors wish to thank two referees, an Associate Editor and Co-Editors Jianqing Fan and John Marden for their constructive comments and useful suggestions. They are grateful to Hengjian Cui, James Hannan, Hira Koul, Regina Liu, Hannu Oja and Ronald Randles for useful discussions.

## REFERENCES

[1] ALEXANDER, K. S. (1984). Probability inequalities for empirical processes and a law of the iterated logarithm. *Ann. Probab.* **12** 1041–1067. MR0757769

[2] BROWN, B. M. and HETTMANSPERGER, T. P. (1987). Affine invariant rank methods in the bivariate location model. *J. Roy. Statist. Soc. Ser. B* **49** 301–310. MR0928938

[3] CHOI, K. and MARDEN, J. (1997). An approach to multivariate rank tests in multivariate analysis of variance. *J. Amer. Statist. Assoc.* **92** 1581–1590. MR1615267

[4] DONOHO, D. L. (1982). Breakdown properties of multivariate location estimators. Ph.D. qualifying paper, Dept. Statistics, Harvard Univ.

[5] DONOHO, D. L. and GASKO, M. (1992). Breakdown properties of multivariate location parameters and dispersion matrices. *Ann. Statist.* **20** 1803–1827. MR1193313

[6] DUDLEY, R. M. (1999). *Uniform Central Limit Theorems.* Cambridge Univ. Press. MR1720712

[7] HETTMANSPERGER, T. P., MÖTTÖNEN, J., and OJA, H. (1998). Affine invariant multivariate rank test for several samples. *Statist. Sinica* **8** 765–800. MR1651508

[8] HUBER, P. J. (1981). *Robust Statistics.* Wiley, New York. MR0606374

[9] LIU, R. Y. (1990). On a notion of data depth based on random simplices. *Ann. Statist.* **18** 405–414. MR1041400

[10] LIU, R. Y. (1992). Data depth and multivariate rank tests. In $L_1$-*Statistical Analysis and Related Methods* (Y. Dodge, ed.) 279–294. North-Holland, Amsterdam. MR1214839

[11] LIU, R. Y. (1995). Control charts for multivariate processes. *J. Amer. Statist. Assoc.* **90** 1380–1387. MR1379481

[12] LIU, R. Y. and SINGH, K. (1993). A quality index based on data depth and multivariate rank tests. *J. Amer. Statist. Assoc.* **88** 252–260. MR1212489

[13] MASSART, P. (1983). Vitesses de convergence dans le théorème central limite pour des processus empiriques. *C. R. Acad. Sci. Paris Sér. I Math.* **296** 937–940. MR0719281

[14] MASSART, P. (1986). Rates of convergence in the central limit theorem for empirical processes. *Ann. Inst. H. Poincaré Probab.-Statist.* **22** 381–423. MR0871904

[15] MASSÉ, J. C. (1999). Asymptotics for the Tukey depth process, with an application to a multivariate trimmed mean. *Bernoulli* **10** 1–23. MR2061438

[16] MÖTTÖNEN, J., HETTMANSPERGER, T. P., OJA, H., and TIENARI, J. (1998). On the efficiency of affine invariant multivariate rank tests. *J. Multivariate Anal.* **66** 118–132. MR1648529

[17] PURI, M. L. and SEN, P. K. (1971). *Nonparametric Methods in Multivariate Analysis.* Wiley, New York. MR0298844

[18] POLLARD, D. (1990). *Empirical Processes: Theorem and Applications.* Institute of Mathematical Statistics, Hayward, California.

[19] RANDLES, R. H. and PETERS, D. (1990). Multivariate rank tests for the two-sample location problem. *Comm. Statist. Theory Methods* **19** 4225–4238. MR1103009

[20] ROUSSON, V. (2002). On distribution-free tests for the multivariate two-sample location-scale model. *J. Multivariate Anal.* **80** 43–57. MR1889832

[21] SEBER, G. A. F. (1977). *Linear Regression Analysis.* Wiley, New York. MR0436482

[22] STAHEL, W. A. (1981). Breakdown of covariance estimators. Research Report 31, Fachgruppe für Statistik, ETH, Zürich.

[23] TOPCHII, A., TYURIN, Y., and OJA, H. (2003). Inference based on the affine invariant multivariate Mann–Whitney–Wilcoxon statistic. *J. Nonparametr. Statist.* **14** 403–414. MR2017477

[24] TUKEY, J. W. (1975). Mathematics and picturing data. In *Proc. Intern. Congr. Math. Vancouver 1974* **2** 523–531. MR0426989

[25] ZUO, Y. (2003). Projection based depth functions and associated medians. *Ann. Statist.* **31** 1460–1490. MR2012822

[26] ZUO, Y. (2004). Projection based affine equivariant multivariate location estimators with the best possible finite sample breakdown point. *Statist. Sinica* **14** 1199–1208. MR2126348

[27] ZUO, Y. and SERFLING, R. (2000). General notions of statistical depth function. *Ann. Statist.* **28** 461–482. MR1790005

DEPARTMENT OF STATISTICS AND PROBABILITY
MICHIGAN STATE UNIVERSITY
EAST LANSING, MICHIGAN 48824
USA
E-MAIL: zuo@msu.edu

DEPARTMENT OF STATISTICS
UNIVERSITY OF ILLINOIS
CHAMPAIGN, ILLINOIS 61820
USA
E-MAIL: x-he@stat.uiuc.edu